\documentclass[11pt]{amsart}

\usepackage{amsmath,amssymb,amsthm,mathtools}
\usepackage{microtype}
\usepackage{xcolor}
\definecolor{mcitecolor}{HTML}{5590B4}
\usepackage[hidelinks]{hyperref}
\hypersetup{
linkcolor  = mcitecolor!85!black,
citecolor  = mcitecolor!85!black,
urlcolor   = blue!85!black,
colorlinks = true,
}
\usepackage{comment}
\usepackage[a4paper,margin=3cm]{geometry}
\usepackage{amsaddr}

\allowdisplaybreaks

\newtheorem{theorem}{Theorem}[section]
\newtheorem{lemma}[theorem]{Lemma}
\newtheorem{proposition}[theorem]{Proposition}
\newtheorem{corollary}[theorem]{Corollary}
\theoremstyle{definition}
\newtheorem{definition}[theorem]{Definition}
\theoremstyle{remark}
\newtheorem{remark}[theorem]{Remark}

\newcommand{\Ftwo}{\mathbb{F}_2}
\newcommand{\Rcal}{\mathcal{R}}
\newcommand{\Qcal}{\mathcal{Q}}
\newcommand{\dd}{\mathbin{\dot\cup}}

\title[A Helly-Type Theorem]
{A Helly-Type Theorem for two-component convex sets}

\author{Giuliamaria Menara}
\address[GM]{\vspace{-3mm}\small Institut de Mathématiques de Jussieu - Paris Rive Gauche}
\email{\href{mailto:giuliamaria.menara@imj-prg.fr} {\nolinkurl{giuliamaria.menara@imj-prg.fr}}}

\begin{document}

\begin{abstract}
In any fixed dimension, consider a finite collection of sets, each consisting of exactly two disjoint, closed, convex pieces. We show that to guarantee the intersection of the whole collection also consists of exactly two such convex pieces, it suffices to verify the same two-piece structure for all intersections of some subfamilies of intermediate size, thus answering a question of Gil Kalai.
\end{abstract}

\maketitle

\section{Introduction}

Helly's theorem \cite[Theorem 1.3.2]{matouvsek2002lectures} is one of the basic local-to-global principles of discrete geometry, and it concerns the intersection of convex sets.
In its classical form, it says that a finite family of convex sets in $\mathbb{R}^d$ has nonempty intersection as soon as every $d+1$ of its members have nonempty intersection. The question considered here asks what remains true when the members of the family are allowed to be disconnected, but in a ``structured'' way.

More precisely, a problem by Gil Kalai in the open problems collection of the 2020 Discrete Geometry meeting~\cite[Problem~13]{adiprasito2021discrete} asks the following. Suppose that every member of a family in $\mathbb{R}^d$ is the disjoint union of two nonempty closed convex sets, and that the
intersection of any $2,3,\dots,d+1$ members has the same property. Must the whole family have nonempty intersection? Closely related formulations appear in~\cite{barany2022helly,KalaiMO}.

Helly-type results for families with a bounded number of components go
back to Grünbaum and Motzkin~\cite{grunbaum1961components}; see also
Larman~\cite{larman1968helly} for a geometric treatment of unions of convex
sets. Amenta~\cite{amenta1996short} proved that if every nonempty finite
intersection is a disjoint union of at most $k$ closed convex sets, then the Helly number is at most $k(d+1)$. Related bounds under
local assumptions on intersections of sufficiently small subfamilies were
obtained by Matoušek~\cite{matouvsek1997helly}. For the abstract combinatorial
principle behind the Grünbaum--Motzkin--Amenta theorem and its history,
see Eckhoff and Nischke~\cite{eckhoff2009morris}.

In this note we answer positively to (a stronger version of) Kalai's question, and the idea of the proof we provide is as follows.
Once the two components of each $F_i$ are labeled by $0$ and $1$, every pair $F_i,F_j$
selects two admissible points of $\Ftwo^2$. Any two-point subset of $\Ftwo^2$ is an affine line, hence it is described by one binary equation.
Therefore, the hypothesis that every triple intersection $F_i\cap F_j \cap F_k$ has exactly two connected components and both components are closed and convex, says that these equations are compatible on every triangle of the complete graph.
At this point, the important observation is that this local condition
leaves exactly one global binary degree of freedom, and this part of the argument reduces to checking the compatibility of local constraints on a labeled simplicial complex.
The result is the following.

\medskip
\noindent\textbf{Main Theorem (Theorem~\ref{thm:main}).}
\emph{Let $d\geq2$, let $n\geq d+2$, and let
$F_1,\ldots,F_n\subseteq\mathbb{R}^d$. Suppose that
$\bigcap_{i\in I}F_i$ has exactly two nonempty closed convex components for
every $I\subseteq[n]$ with $|I|\in\{1,2,3,d+1\}$.
Then $\bigcap_{i=1}^nF_i$ has exactly two nonempty closed convex components.}

\medskip
Thus the conclusion is stronger than the nonemptiness asked for in the problem, and the assumptions on intersections of $4,\dots,d$ members are not used.
The restriction $d\geq2$ is real: in dimension one the local condition of the original problem only concerns pairs, and a four-set counterexample is given in Proposition~\ref{prop:d1}.

The note is organized as follows. Section~\ref{sec:setup} records the
component decomposition that allows us to pass from geometry to binary
labels. Section~\ref{sec:binary} contains the local-to-global principle for the
binary constraints. We prove the main theorem in Section~\ref{sec:proof},
and conclude with the one-dimensional counterexample and a compactness
remark in Section~\ref{sec:remarks}.

\section{Setup}
\label{sec:setup}

We begin with the terminology used throughout the note.

\begin{definition}\label{def:two-component}
A subset $F\subseteq\mathbb{R}^d$ has \emph{exactly two closed convex
components} if it has exactly two connected components and both of them are
nonempty closed convex sets.
\end{definition}

After arbitrarily labeling the two components, we write $F=C^0\dd C^1$, where the labels $0$ and $1$ will only be used to encode choices by elements of
$\Ftwo$.
The following Lemma \ref{lem:cells} shows that this encoding identifies the connected components of a finite intersection with its nonempty binary cells.

\begin{lemma}
\label{lem:cells}
Let $I$ be a nonempty finite set. For every $i\in I$, let $F_i=C_i^0\dd C_i^1$, where $C_i^0$ and $C_i^1$ are disjoint closed convex sets. For
$x\in\Ftwo^I$, set
\[
 D_x:=\bigcap_{i\in I}C_i^{x_i}.
\]
Then the nonempty sets $D_x$ are precisely the connected components of
$\bigcap_{i\in I}F_i$.
\end{lemma}

\begin{proof}
Distributivity gives a finite disjoint decomposition
\begin{equation}
\label{eq:cell-decomposition}
\bigcap_{i\in I}F_i
=
\mathop{\dd}_{\substack{x\in\Ftwo^I\\D_x\neq\emptyset}}D_x.
\end{equation}
Every nonempty $D_x$ is closed and convex, hence connected. Two different
cells are disjoint, because two different assignments disagree at some
index $i$ and $C_i^0\cap C_i^1=\emptyset$.

It remains to see that two cells cannot belong to the same connected
component of the union. Since the decomposition in
\eqref{eq:cell-decomposition} is finite and every cell is closed, the complement of any one cell is closed in the union.
Hence every cell is both open and closed in the relative topology.
A connected subset of the union can therefore meet at most one cell, and this proves the claim.
\end{proof}

We record the resulting binary data for later use. Given a family
$F_i=C_i^0\dd C_i^1$, and a set $I$ of indices, define
\begin{equation}\label{eq:R-definition}
 \Rcal(I):=
 \left\{
 x\in\Ftwo^I:
 \bigcap_{i\in I}C_i^{x_i}\neq\emptyset
 \right\}.
\end{equation}
By Lemma~\ref{lem:cells}, $|\Rcal(I)|$ is exactly the number of connected
components of $\bigcap_{i\in I}F_i$.

\section{The binary local-to-global principle}
\label{sec:binary}

We now forget the geometry for a moment and focus only on the pairwise constraints.
The key fact here is that the local intersection hypotheses force a global system of binary equations to have exactly two solutions. This section will be mostly linear algebra over $\Ftwo$.  

Let $V$ be a finite set. For every pair $\{i,j\}\subseteq V$, let $L_{ij}\subseteq\Ftwo^{\{i,j\}}$ be a two-element set. For $I\subseteq V$, define
\begin{equation}\label{eq:Q}
\mathcal{Q}(I):=
\left\{
 x\in\Ftwo^I:
 x|_{\{i,j\}}\in L_{ij}
 \text{ for every }\{i,j\}\subseteq I
\right\}.
\end{equation}
Every two-element subset of $\Ftwo^2$ is an affine line. Thus there is a unique equation, of one of the three forms
\[
x_i=c,
\qquad
x_j=c,
\qquad\text{or}\qquad
x_i+x_j=c,
\qquad c\in\Ftwo,
\]
whose solution set is $L_{ij}$.

\begin{definition}
Call an edge $ij$ a \emph{parity edge} if its equation is $x_i+x_j=c_{ij}$.
\end{definition}

Every other edge fixes one of its endpoints. We orient such an edge toward the endpoint that it fixes. Thus, for example, $i\to j$ means that the equation on $ij$ is $x_j=c_{ij}$.

We will need three technical results for the local-to-global principle: Lemma~\ref{lem:parity-classes} identifies the parity classes, Lemma~\ref{lem:class-tournament} orders them transitively, and finally Lemma~\ref{lem:two-global-assignments} yields exactly two global assignments.

\begin{remark}
The idea of the following Lemma \ref{lem:parity-classes} is that, even though a system of pairwise binary constraints can fail globally due to an odd cycle, e.g. $x_1=x_2$, $x_2=x_3$, $x_1 \neq x_3$, in our case the hypothesis on triples rules out such an obstruction, because on every triangle the three equations form a consistent affine system of rank two.
\end{remark}

\begin{lemma}
\label{lem:parity-classes}
Let $V$ be finite, with $|V|\geq 3$, and suppose that $|\mathcal{Q}(T)|=2$ for every three-element subset $T\subseteq V$. Then the graph formed by the parity edges is a disjoint union of complete graphs.
\end{lemma}

\begin{proof}
We have to show that the relation $i\sim j$ defined by $i=j$ or by $ij$ being a parity edge is an equivalence relation on $V$.

For every triple, the corresponding affine system is consistent and has exactly two solutions. Its coefficient matrix over $\Ftwo$ therefore has rank exactly two.

A triangle cannot contain exactly two parity edges. In fact, suppose that $ij$ and $ik$ are parity edges while $jk$ is not. The coefficient vectors of the three equations are then
\[
\{e_i+e_j,\quad e_i+e_k,\quad e_j\}
\quad \text{or} \quad 
\{e_i+e_j,\quad e_i+e_k,\quad e_k\}.
\]
In either case they are linearly independent, contradicting rank two. Consequently, if $ij$ and $ik$ are parity edges, then $jk$ is a parity edge as well. Thus $\sim$ is transitive; reflexivity and symmetry are immediate. Hence parity adjacency is an equivalence relation, and each equivalence class spans a complete graph of parity edges.
\end{proof}

Let $V_1,\ldots,V_m$ denote the parity classes given by Lemma~\ref{lem:parity-classes}. All edges inside a class are parity edges, and all edges between different classes are oriented edges.

\begin{lemma}
\label{lem:class-tournament}
Assume the hypotheses of Lemma~\ref{lem:parity-classes}. Between any two parity classes, all edges have one common direction, and the induced tournament on the classes is transitive.
\end{lemma}

\begin{proof}
First, all edges between two fixed classes $A$ and $B$ have the same class-level direction. To see this, fix $a,a'\in A$ and $b\in B$. In the triangle $a,a',b$, the edge $aa'$ is a parity edge. If one of $ab,a'b$ fixed its endpoint in $A$ and the other fixed $b$, then the three coefficient vectors would have rank three. Thus, for fixed $b$, the direction does not depend on the choice of the vertex in $A$. The analogous argument in a triangle $a,b,b'$, with $b,b'\in B$, shows that it does not depend on the choice of the vertex in $B$. We may therefore write $A\to B$ when every edge between $A$ and $B$ fixes its endpoint in $B$.

This defines a tournament on the parity classes. It has no directed triangle. Indeed, if $A\to B$, $B\to C$, $C\to A$, and $a\in A$, $b\in B$, $c\in C$, then the equations on $ab,bc,ca$ fix $x_b,x_c,x_a$, respectively. Their coefficient vectors are $e_b,e_c,e_a$, of rank three, contradicting the hypothesis on the triple $\{a,b,c\}$. Hence the tournament is transitive. After renumbering the classes, we may assume that
\begin{equation}\label{eq:class-order}
V_1\to V_2\to\cdots\to V_m,
\end{equation}
meaning that, when $r<s$, every edge between $V_r$ and $V_s$ fixes the variable in $V_s$.
\end{proof}

With the parity classes linearly ordered, the problem becomes straightforward: the first class has a free binary choice (and this will eventually correspond to the two connected components of the full intersection), and every variable in the later classes is then forced by the directed edges from earlier classes.

\begin{lemma}
\label{lem:two-global-assignments}
Under the hypotheses of Lemma~\ref{lem:parity-classes}, it holds that $|\mathcal{Q}(V)|=2$.
\end{lemma}

\begin{proof}
We keep the notation from Lemmas~\ref{lem:parity-classes} and~\ref{lem:class-tournament}: the parity classes are $V_1,\dots,V_m$, ordered so that, whenever $r<s$, every edge between $V_r$ and $V_s$ fixes the variable in $V_s$.

We first analyze each parity class independently.

Fix a parity class $V_r$. For each pair of distinct vertices $u,v\in V_r$, the corresponding equation has the form $x_u+x_v=c_{uv}$.
For any three distinct vertices $u,v,w\in V_r$, the system on $\{u,v,w\}$ is consistent. Summing its three equations gives
\begin{equation}
\label{eq:coboundary}
c_{uv}+c_{vw}+c_{uw}=0.
\end{equation}
Thus the labels $c_{uv}$ form a $1$-cocycle on the simplex with vertex set $V_r$. Since a simplex has trivial first cohomology, this cocycle is a coboundary, meaning there exists a function
\[
\tau_r\colon V_r\longrightarrow\mathbb F_2
\]
such that $c_{uv}=\tau_r(u)+\tau_r(v)$ for all distinct $u,v\in V_r$.
Concretely, one may choose a base vertex $a_r\in V_r$, set
$\tau_r(a_r)=0$, and define
\[
\tau_r(v)=c_{a_rv}
\qquad (v\neq a_r);
\]
then equation~\ref{eq:coboundary} gives the desired identity. The same conclusion is immediate when $|V_r|\leq 2$.
It follows that the equations inside $V_r$ are equivalent to requiring
that $x_v+\tau_r(v)$ are constant on $V_r$. Hence their solutions are
exactly
\begin{equation}
\label{eq:solutions}
x_v=t_r+\tau_r(v),
\qquad v\in V_r,\quad t_r\in\mathbb F_2.
\end{equation}
Thus each parity class carries one binary parameter before the constraints
between different classes are imposed.

We now show that only the first parameter remains free.

If $m=1$, equation~\ref{eq:solutions} already gives exactly two global assignments.
Assume therefore that $m\geq 2$.

Let $s\geq 2$ and $v\in V_s$. For every $u\in V_1$, the edge
$uv$ fixes $x_v$. If $u,u'\in V_1$ are distinct, consistency of
the system on the triangle $\{u,u',v\}$ implies that the edges $uv$
and $u'v$ prescribe the same value for $x_v$. When $V_1$ is a
singleton, this is immediate. Denote this common value by $b_v$.
We claim that the assignments
\begin{equation}
\label{eq:assignments}
x_v=
\begin{cases}
t+\tau_1(v), & v\in V_1,\\
b_v, & v\in V_2\cup\cdots\cup V_m,
\end{cases}
\qquad t\in\mathbb F_2,
\end{equation}
satisfy all pair constraints.

We already showed that the constraints inside $V_1$ are satisfied, and the constraints between $V_1$ and every later class are satisfied by the definition of $b_v$. It remains to check the constraints whose two endpoints lie in later classes.

First, let $v,w\in V_s$ be distinct, with $s\geq2$, and choose
$u\in V_1$. On the triangle $\{u,v,w\}$, the edges $uv$ and $uw$ fix $x_v=b_v$  and $x_w=b_w$, while $vw$ is a parity edge. Since this triangle system is consistent, $b_v+b_w=c_{vw}$.
Thus all parity equations inside $V_s$ are satisfied.

Next, let $1<r<s$, let $w\in V_r$, and let $v\in V_s$. Choose
$u\in V_1$. On the triangle $\{u,w,v\}$, both $uv$ and $wv$
fix $x_v$. The first implies $x_v=b_v$, and consistency forces the second to prescribe the same value. Hence every constraint between $V_r$ and $V_s$ is also satisfied.

Therefore each of the two values $t\in\mathbb F_2$ in equation \ref{eq:assignments} gives a
global solution. Conversely, every global solution satisfies the equations
inside $V_1$, so its restriction to $V_1$ is determined by a unique
choice of $t$, and every variable outside $V_1$ is then forced to equal $b_v$. Thus there are no other global solutions, meaning $\mathcal Q(V)|=2$.
\end{proof}

\begin{proposition}[Binary local-to-global principle]\label{prop:binary-local-global}
Let $V$ be finite, with $|V|\geq 3$. Suppose that $|\mathcal{Q}(T)|=2$ for every three-element subset $T\subseteq V$. Then $|\mathcal{Q}(I)|=2$ for every $I\subseteq V$ with $|I|\geq 2$.
\end{proposition}

\begin{proof}
For $|I|=2$, the conclusion is part of the definition of $L_{ij}$. If $|I|\geq 3$, the hypothesis is inherited by the subsystem induced on $I$. Applying Lemmas~\ref{lem:parity-classes}, \ref{lem:class-tournament}, and \ref{lem:two-global-assignments} to that subsystem gives $|\mathcal{Q}(I)|=2$.
\end{proof}

\section{Proof of the main theorem}
\label{sec:proof}

We return to the family $F_1,\ldots,F_n\subseteq\mathbb{R}^d$. Label the two components of $F_i$ as $F_i=C_i^0\dd C_i^1$, and let $\Rcal(I)$ be as in~\eqref{eq:R-definition}. For every pair $\{i,j\}$, set $L_{ij}:=\Rcal(\{i,j\})$.

By the hypotheses and Lemma~\ref{lem:cells}, $L_{ij}$ has exactly two
elements. Let $\Qcal(I)$ be the set of assignments satisfying all these
pair constraints, as in~\eqref{eq:Q}. By construction,
\begin{equation}\label{eq:R-subset-Q}
 \Rcal(I)\subseteq\Qcal(I)
 \qquad\text{for every }I\subseteq[n].
\end{equation}

The following short observation is what allows us to apply Proposition \ref{prop:binary-local-global}.
The idea is that we must check that the geometric assumption on triple intersections matches the algebraic assumption of Proposition \ref{prop:binary-local-global}. 

\begin{lemma}\label{lem:triple-equality}
For every three-element set $T\subseteq[n]$, one has $\Rcal(T)=\Qcal(T)$, and hence $|\Qcal(T)|=2$.
\end{lemma}

\begin{proof}
The three pair equations on $T=\{i,j,k\}$ have nonzero coefficient vectors
supported, respectively, in $\{i,j\}$, $\{i,k\}$, $\{j,k\}$.
Their rank cannot be one: over $\Ftwo$, three nonzero vectors of rank one
would be equal, but no nonzero vector is supported in the intersection of
these three sets. Thus the rank is at least two.

On the other hand, the hypothesis and Lemma~\ref{lem:cells} give
$|\Rcal(T)|=2$. Because of~\eqref{eq:R-subset-Q}, the affine system defining
$\Qcal(T)$ is consistent and has at least two solutions. An affine system
in three binary variables whose coefficient rank is at least two has at most
two solutions. Therefore $|\Qcal(T)|=2$, and the inclusion
$\Rcal(T)\subseteq\Qcal(T)$ is an equality.
\end{proof}

\begin{theorem}\label{thm:main}
Let $d\geq2$, let $n\geq d+2$, and let
$F_1,\ldots,F_n\subseteq\mathbb{R}^d$. Suppose that
$\bigcap_{i\in I}F_i$ has exactly two nonempty closed convex components for
every $I\subseteq[n]$ with $|I|\in\{1,2,3,d+1\}$.
Then $\bigcap_{i=1}^nF_i$ has exactly two nonempty closed convex components.
\end{theorem}

\begin{proof}
By Lemma~\ref{lem:triple-equality}, the pair system satisfies the hypotheses
of Proposition~\ref{prop:binary-local-global}. Hence
\begin{equation}\label{eq:Q-two}
 |\Qcal(I)|=2
 \qquad\text{for every }I\subseteq[n]\text{ with }|I|\geq2.
\end{equation}
In particular, $\Qcal([n])$ consists of exactly two global assignments.

Fix $x\in\Qcal([n])$, and let $I\subseteq[n]$ have cardinality $d+1$. Then
$x|_I\in\Qcal(I)$. By hypothesis and Lemma~\ref{lem:cells},
$|\Rcal(I)|=2$, while~\eqref{eq:Q-two} gives $|\Qcal(I)|=2$. Together with
\eqref{eq:R-subset-Q}, this implies $\Rcal(I)=\Qcal(I)$.
Consequently,
\[
 \bigcap_{i\in I}C_i^{x_i}\neq\emptyset
 \qquad\text{for every }I\in\binom{[n]}{d+1}.
\]
The sets $C_i^{x_i}$ are closed and convex. Helly's theorem \cite[Theorem 1.3.2]{matouvsek2002lectures} therefore gives $\bigcap_{i=1}^n C_i^{x_i}\neq\emptyset$. Thus $x\in\Rcal([n])$.

We have proved $\Qcal([n])\subseteq\Rcal([n])$. The reverse inclusion is
\eqref{eq:R-subset-Q}, so $\Rcal([n])=\Qcal([n])$.
The right-hand side has exactly two elements, and Lemma~\ref{lem:cells} shows that $\bigcap_{i=1}^nF_i$ has exactly two nonempty closed convex components.
\end{proof}

\begin{corollary}[{\cite[Problem~13]{adiprasito2021discrete}}]
\label{cor:problem13}
Let $d\geq2$, and let $\mathcal{F}$ be a finite family of subsets of
$\mathbb{R}^d$. Suppose that every member of $\mathcal{F}$ is the disjoint
union of exactly two nonempty closed convex sets, and that the same holds for
the intersection of every $k$ members, for each $2\leq k\leq d+1$. Then the
intersection of all members of $\mathcal{F}$ has exactly two nonempty closed
convex components. In particular, it is nonempty.
\end{corollary}


\section{The one-dimensional case and compactness}
\label{sec:remarks}

The triple condition in Proposition~\ref{prop:binary-local-global} cannot be omitted. This is visible in dimension one, where $d+1=2$ and the original local
condition only concerns pairs.

\begin{proposition}
\label{prop:d1}
There are four subsets of $\mathbb{R}$, each having exactly two nonempty
closed convex components, such that every pairwise intersection has exactly
two nonempty closed convex components, while the intersection of all four
sets is empty.
\end{proposition}

\begin{proof}
Let
\[
\begin{aligned}
 F_1=[0,1]\dd[2,4],\quad
 F_2=[0,1]\dd[3,5],\quad
 F_3=[0,2]\dd[5,6],\quad
 F_4=[2,3]\dd[4,5].
\end{aligned}
\]
The six pairwise intersections are
\[
\begin{array}{lll}
 F_1\cap F_2=[0,1]\dd[3,4],
 &
 F_1\cap F_3=[0,1]\dd\{2\},
 &
 F_1\cap F_4=[2,3]\dd\{4\},
 \\[1mm]
 F_2\cap F_3=[0,1]\dd\{5\},
 &
 F_2\cap F_4=\{3\}\dd[4,5],
 &
 F_3\cap F_4=\{2\}\dd\{5\}.
\end{array}
\]
Thus every pairwise intersection has exactly two closed convex components.
However, $F_1\cap F_2\cap F_3=[0,1]$ and $[0,1]\cap F_4=\emptyset$, so $F_1\cap F_2\cap F_3\cap F_4=\emptyset$.
\end{proof}

\begin{remark}[Infinite families]\label{rem:infinite}
Closedness alone does not allow one to pass from finite to infinite
families. In $\mathbb{R}^2$, let 
\[
F_m=
 \bigl([m,\infty)\times[1,\infty)\bigr)
 \dd
 \bigl([m,\infty)\times(-\infty,-1]\bigr),
 \qquad m\geq1.
\]
Every nonempty finite intersection again has exactly two nonempty closed
convex components, while $\bigcap_{m\geq1}F_m=\emptyset$.

If all components are compact, the usual finite-intersection argument applies.
\end{remark}

\subsection*{Acknowledgements}
\sloppy
I am grateful to Karim Adiprasito for valuable comments, and to Zuzana Pat\'akov\'a for thoughtful feedback on this manuscript.\\
The author was supported by Horizon Europe ERC Grant number: 101045750 / Project acronym: HodgeGeoComb.

\bibliographystyle{alphaurl}
\bibliography{bibliography}

\end{document}